\documentclass[12pt,letterpaper]{article}
\usepackage[utf8]{inputenc}
\usepackage[english]{babel}
\usepackage{amsmath}
\usepackage{amsfonts}
\usepackage{amssymb}
\usepackage{makeidx}
\usepackage{graphicx}
\usepackage{amsthm}
\usepackage[none]{hyphenat}
\usepackage[colorlinks=true,citecolor=red,linkcolor=blue,pdfpagetransition=Blinds]{hyperref}
\theoremstyle{plain}
\newtheorem{theorem}{Theorem}[section]

\newtheorem{lemma}[theorem]{Lemma}

\theoremstyle{definition}
\newtheorem*{definition}{Definition}

\theoremstyle{remark}

\def\r{\mathbb{R}}
\def\rn{\mathbb{R}^N}

\def\o{\Omega}
\def\io{\int_{\Omega}}
\def\co{\overline{\Omega}}
\usepackage[left=2.54cm,right=2.54cm,top=2.54cm,bottom=2.54cm]{geometry}
\title{Sobolev embeddings on domains involving two types of symmetries}
\author{
David Flores-Flores$^{1}$\\
Alfredo Cano$^{2}$\\
Eric Hernández-Martínez$^{3}$.\\ \\
\small{\textit{$^{1}$ Instituto de Matemáticas, Universidad Nacional Autónoma de México, Circuito Exterior,}}\\
\small{\textit{Ciudad Universitaria, 04510, Coyoacán, Ciudad de México, Mexico}}\\
\small{\textit{davhfl@hotmail.com}}\\
\small{\textit{$^{2}$Facultad de Ciencias, Universidad Autónoma del Estado de México,}}\\
\small{\textit{Campus El Cerrillo Piedras Blancas, Toluca, 50200, Estado de México}}\\
\small{\textit{calfredo420@gmail.com}}\\ 
\small{\textit{$^{3}$ Colegio de Ciencia y Tecnología, Academia de Matemáticas,}}\\
\small{\textit{Universidad Autónoma de la Ciudad de México,}}\\
\small{\textit{Calle Prolongación San Isidro No. 151, Col. San Lorenzo Tezonco,}}\\
\small{\textit{Del. Iztapalapa, 09790, Ciudad de México}}\\
\small{\textit{ebric2001@hotmail.com}}\\
}
\date{}
\begin{document}
\maketitle
\begin{abstract}
It is well known that Sobolev embeddings can be improved in the presence of symmetries. In this article, we considere the situation in which given a domain $\Omega=\Omega_1 \times \Omega_2$ in $\mathbb{R}^N$ with a cylindrical symmetry, and acting a group $G$ in $\Omega_1$, for this situation it is shown that the critical Sobolev exponent increases in the case of embeddings into weighted spaces $L^{q}_{h}(\Omega)$. In this paper, we will enunciate several results related to compact embeddings of a Sobolev space with radially symmetric functions into some weighted space $L^{q}$, with $q$ higher than the usual critical exponent.
\end{abstract}
\section{Introduction}\label{sec1}
In the field of partial differential equations, the study of Sobolev embeddings takes a vital role for the variational method. The symmetries also gives important properties in variational problems, so the study of Sobolev embeddings with symmetries, allowing to solve more problems related to partial differential equations. It can be seen in recently results, for instance  (\cite{ackermann2017spectral}, \cite{edmunds2000sobolev}, \cite{edmunds2002sobolev}, \cite{lions1982symetrie}, \cite{musina2014weighted}). It is well known that Sobolev spaces can be improved if we deal with subspaces of invariant functions under the action of some group of symmetries (see \cite{aubin2012nonlinear}, \cite{cho2009sobolev}, \cite{faget2002best}, \cite{fan2001compact}, \cite{faraci2011low}, \cite{skrzypczak2003heat}, \cite{willem1997minimax}). These conditions have been used to prove the existence of solutions to value problems at the boundary. In some cases, it is sufficient to consider invariant functions to analyze under what conditions the compactness of these embeddings occurs. However, symmetry alone does not provide the desired improvements if we restrict ourselves to the usual $L^p$ spaces.\\
\\Given a domain $\Omega$ in $\mathbb{R}^N$, is well known that for $1<p<\infty$, if $p<N$, the Sobolev space $W^{1,p}(\Omega)$ can be embedded continuously into $L^q(\Omega)$ if $p<q \leq p^{*}$, where $p^{*}=\frac{Np}{n-p}$, (see  \cite{adams2003sobolev}, \cite{brezis2011functional}) (remembering that if $p=2$, this spaces are Hilbert spaces and are usually denoted as $H^{1}$). Also, if $\Omega$ is bounded and $q<p*$, the embedding is compact. If we take $W^{1,p}_{s}(\Omega)=\{u \in W^{1,p}(\Omega): u \; \mbox{is radially symmetric, that is,} \; u(x)=u(\Vert x \Vert), \; \forall x \in \Omega \}$, by the Lions result (see \cite{lions1982symetrie}, \cite{lions1985concentration}), $W^{1,p}_{s}(\Omega)$ can be embedded compactly into $L^q(\Omega)$ for $p<q<p^{*}$. The idea here is to show how that Sobolev embeddings can be improved in the presence of symmetries by group actions. This includes embeddings in some weighted spaces $L^{p}$, with $p$ higher than the usual critical Sobolev exponent. Such phenomena have been observed in specific contexts by various authors. On the one hand, Hebey and Vaugon (\cite{Hebey1996}, \cite{HEBEY1997859}) studied this situation in the context of a compact Riemannian manifold $M$ and a compact group of isometries $G$. They considered the dimension of the group orbit, denoted by $k$, and they analyzed the Sobolev Space $W^{1,p}_{0,G}(M)= \{ u \in W^{1,p}_{0}(M): \; u \circ \sigma = u  \; \; \forall \; \sigma \in G\}$. They proved that the critical Sobolev exponent increases if there are no points on an open subset with closure compact $\Omega$ of $M$ with discrete orbits under the action of $G$, they obtained a greatest critical exponent for the embedding $W^{1,p}_{0,G}(\Omega) \hookrightarrow L^{q}(\Omega)$ for $1<q<\frac{p(N-k)}{N-k-p}$.\\
\\On the other hand, Wang (\cite{wang2006sobolev}) established embedding results in a regular domain with cylindrical symmetry. Making $H^1_s(\Omega) = \{ u \in H^1_0(\Omega) : u(x_1, x_2)=w(x_1, \Vert x_2 \Vert ), \; \forall x_2 \in \Omega_2 \}$, he proved that taking $h=\Vert x \Vert ^l$, $x \in \Omega_2$ and $l$ a real positive number, which is radially symmetric on $\Omega_2$, then, there exists a positive number $\tau$, that depends of dimensions of $\Omega_1$ and $\Omega_2$ such that the embedding $H^1_s(\Omega) \hookrightarrow L^{q}_{h}(\Omega)$ is compact for $1 < q < 2^{*} + \tau$, where $2^{*}$ is the one given by the Sobolev embedding theorems.\\
\\Our goal here is the study of Sobolev embeddings, generalizing Wang's results using those obtained by Hebey and Vaugon. First, we treat Sobolev embeddings in a domain with cylindrical symmetry. This domain will be taken as follows: let $\Omega = \Omega_1 \times \Omega_2$ in $\mathbb{R}^N$; where $\Omega_1$ is a bounded set of class $C^1$ in $\mathbb{R}^{m}$, and $\Omega_2$ in $\mathbb{R}^{j}$ is a ball of radius $R$ centered at the origin. After this, we will involve the action of the groups, making a group $G$ acts on $\Omega_1$. With this, we define a $G-$invariant on $\Omega_1$ and radially symmetric on $\Omega_2$ Sobolev space $H^1_{s,G} = \{ u \in H^1_s(\Omega) : \forall g \in G, u(gx_1, x_2) =u(x_1,x_2) \}$, and next, to $y \in \Omega_2$, we take a continuous function $h(y)$, and consider the weighted Lebesgue space $L^{q}_{h}(\Omega)$. Based on this, we will prove that under certain conditions $H^{1}_{s,G}(\Omega)$ can be embedded compactly into $L^{q}_{h}(\Omega)$ for $1 < q < 2^{*}_{k} + \psi$ where $\psi$ is greater than the term $\tau$ given by Wang and $2^{*}_{k}=\frac{2(n-k)}{n-k-2}$ if $n-k>2$. \\
\\Our first result is 
\begin{theorem} \label{teo1}
Let $\Omega = \Omega_1 \times \Omega_2$ with dim $\Omega=N$, dim $\Omega_1 =m \geq 2$ and dim $\Omega_2=j \geq 3$, $G$ a compact subgroup of $\mathcal{O}(m)$ acting on $\Omega_1$ and $k = \displaystyle \min_{x \in \overline{\Omega_1}} \mbox{dim} \; O^{x}_{G}$ with $k \geq 1$. Suppose that $h=\Vert x_2 \Vert^l$, with $l >0$ a real number. Then the embedding $H^1_{s,G}(\Omega) \hookrightarrow L^q_h(\Omega)$ is compact for $q \in (1, 2^{*}_{k} + \psi)$, where $\psi = \frac{2}{N-k-2} \min \left\lbrace \frac{2(j-2)}{m-k},l \right\rbrace$.
\end{theorem}

In the case of lower dimension for $\Omega_2$ we give the next result. 
\begin{theorem} \label{teo2}
Assume that dim $\Omega_1 =m \geq 2$ and dim $\Omega_2=j = 2$. Let $G$ a compact subgroup of $\mathcal{O}(m)$ acting on $\Omega_1$ and $k = \displaystyle \min_{x \in \overline{\Omega_1}} \mbox{dim} \; O^{x}_{G}$ with $k \geq 1$. Let $h=\Vert x_2 \Vert^l$, with $l >0$ a real number. Then the embedding $H^1_{s,G}(\Omega) \hookrightarrow L^q_h(\Omega)$ is compact for $q \in (1, 2^{*}_k + \psi)$, where
\begin{equation*}
\psi = -\frac{2k(N-k)}{m(m-k)(m-k+1)} +\frac{2(m+1)}{m(m-k+1)} \min \left\lbrace l, \frac{1}{m-k} \right\rbrace.
\end{equation*}
\end{theorem}

This paper consists of four sections. Section \ref{sec2} contains a profile and some previous results obtained by Wang; in Section \ref{sec3} we will treat some preliminaries; the proofs for the main results are given in Section \ref{sec4}.
\section{Some previous results}\label{sec2}
We will introduce some previous results that we will use in this article. We begin with some notations and results that were enunciated by Wang \cite{wang2006sobolev}.
\begin{definition} \label{SpaceofWang}
Let $\Omega= \Omega_1 \times \Omega_2$ in $\mathbb{R}^{N}$; where $\Omega_1$ is a bounded set of class $C^1$ in $\mathbb{R}^{m}$, and $\Omega_2$ in $\mathbb{R}^{j}$ is a ball of radius $R$ centered at the origin. We denote
\begin{equation*}
H^1_s = \{ u \in H^1_0(\Omega) : u(x_1, x_2)=u(x_1, \Vert x_2 \Vert ), \; \forall x_2 \in \Omega_2 \}.
\end{equation*}
\end{definition}
\begin{definition}
Let $h$ be a nonnegative Hölder continuous function on $\overline{\Omega}$, denote by $L^p_h(\Omega)$ the Banach
space, or weighted $L^p$ space, with the norm
\begin{equation*}
\Vert u \Vert_{L^p_h(\Omega)} = \left( \displaystyle \int_{\Omega} h \vert u \vert ^p \right)^{\frac{1}{p}}.
\end{equation*}
\end{definition}
Denote for $m$, $j$ and $N$ the dimensions of $\Omega_1$, $\Omega_2$ and $\Omega$ respectively, then $N=m+j$; and
\begin{equation*}
\begin{split}
\nabla_1 u & = \left(\frac{\partial u}{\partial x_{1_1}}, \frac{\partial u}{\partial x_{1_2}}, \ldots , \frac{\partial u}{\partial x_{1_m}} \right),\\
\nabla_2 u & = \left(\frac{\partial u}{\partial x_{2_1}}, \frac{\partial u}{\partial x_{2_2}}, \ldots , \frac{\partial u}{\partial x_{2_j}} \right)\\
\end{split}
\end{equation*}
are the partial gradient on $\Omega_1$ and $\Omega_2$ respectively, and $\nabla u=(\nabla_1 u, \nabla_2 u)$.\\
\\Next, we will introduce the following lemmas (see Wang \cite{wang2006sobolev}).
\begin{lemma}[Wang, 2006] \label{Lema1W}
Let dim $\Omega_2 =j \geq 3$, then for $u \in  H^{1}_{s}(\Omega)$, we have
\begin{equation*}
\vert u(x) \vert \leq \dfrac{C_j}{\Vert x_2 \Vert^{\frac{j-2}{2}}} \left( \displaystyle \int_{\Omega_2} \vert \nabla  u \vert^2 \right)^{\frac{1}{2}}
\end{equation*}
for all $x \in  \Omega$, where $C_j$ is a constant that depends only of $j$.
\end{lemma}
In the case of $j=2$, we need this lemma
\begin{lemma}[Wang, 2006] \label{Lema2W}
Let dim $\Omega_2=j=2$, and let $b>1$ be real number, then for $u \in H^{1}_{0}(\Omega)$, and for all $x \in \Omega$, we have
\begin{equation} \label{Desig1}
\vert u(x) \vert \leq \dfrac{A_j}{\Vert x_2 \Vert^{j-1}} \displaystyle \int_{\Omega_2} \vert \nabla u  \vert,
\end{equation}
\begin{equation*}
\vert u(x) \vert^b \leq \dfrac{b A_j}{\Vert x_2 \Vert^{j-1}} \displaystyle \int_{\Omega_2} \vert u \vert^{b-1} \vert \nabla u  \vert.
\end{equation*}
\end{lemma}

\section{Background material} \label{sec3}
We present here the notations and the background material that we will use below. We denote by $\mathcal{O}(n)$ the orthogonal group in dimension $n$. A subgroup $G \subset \mathcal{O}(n)$ acts on the $G-$invariant domain $\Omega$ and it induces a well defined action on the function space with 
\begin{equation*}
(gu)(x)=u(g^{-1}x).
\end{equation*}
If $G$ denotes a compact subgroup of $\mathcal{O}(n)$, for a given $G-$invariant open set $\Omega$ of class $C^1$, we set
\begin{equation*}
L^{p}_{G}(\Omega) = \{u \in L^{p}(\Omega) : gu=u, \forall g \in G \}
\end{equation*}
with the usual norm
\begin{equation*}
\Vert u \Vert_{L^{p}_{G}(\Omega)} = \left( \displaystyle \int_{\Omega} \vert u \vert ^p \right)^{\frac{1}{p}}.
\end{equation*}
Let us now recall some facts about the group action.
\begin{definition}
For $G$ a group and $x$ a point of $\Omega$, we denote by
\begin{equation*}
O^{x}_{G}= \{gx: g \in G \}
\end{equation*}
the orbit of $x$ under the action of $G$.
\end{definition}
Taking $\Omega = \Omega_1 \times \Omega_2$ with the conditions described in the definition \ref{SpaceofWang}; and considering $G$ a compact subgroup of $\mathcal{O}(m)$ acting on $\Omega_1$. We denote by
\begin{equation*}
H^1_{s,G} = \{ u \in H^1_s(\Omega) : \forall g \in G, u(gx_1, x_2) =u(x_1,x_2) \}.
\end{equation*}
Next, we give the following definition, also it can be found in \cite{HEBEY1997859}, \cite{ivanov2007weighted}.
\begin{definition}
For $p \geq 1$ a real number, $n$ and $k$ two positive integers, with $n \geq 2$ and $0 \leq k \leq n-1$. We define the critical Sobolev exponent $p_k^{*}= p^{*}(n,k,p)$
\begin{equation*}
p_k^{*}=  \left\lbrace \begin{matrix}
\frac{(n-k)p}{n-k-p} & if & n-k>p\\
\infty & if & n-k \leq p.
\end{matrix} \right.
\end{equation*}
\end{definition}
We remark that if $k=0$, we have the conventional critical exponent $p^{*}=\frac{np}{n-p}$, and $p_k^{*} > p^{*}$.\\
\\We enunciate the first lemma that it helps us in our results, and was proved in \cite{sobolev1938theoreme}.
\begin{lemma}[Sobolev's lemma] \label{LemaSobolev}
Let $p>1$ and $q>1$ two real numbers. Define $\lambda$ by $\frac{1}{p}+ \frac{1}{q}+\frac{\lambda}{N}=2$ with $0<\lambda<N$. Then, there exists a constant $K(p,q,N)$, such that for all $f \in L^p(\Omega)$ and $g \in L^q(\Omega)$,
\begin{equation*}
\displaystyle \int_{\Omega} \displaystyle \int_{\Omega} \frac{
f(x)g(y)}{\Vert x-y \Vert^{\lambda}} \, dx \, dy = \displaystyle \int_{\Omega} g(y) \displaystyle \int_{\Omega} \frac{f(x)}{\Vert x-y \Vert^{\lambda}} \, dx \, dy \leq K(p,q,N) \Vert f \Vert_{L^p(\Omega)} \Vert g \Vert_{L^q(\Omega)}.
\end{equation*}
\end{lemma}
Based in lemma \ref{LemaSobolev}, now we define the following operator and its properties.
\begin{definition}
Let $1 < p < \infty$ and $\lambda$ a real number with $0 < \lambda <N$. Given $u \in L^p(\Omega)$, we define
\begin{equation} \label{V function}
H(u)(y)= \displaystyle \int_{\Omega} \frac{u(x)}{\Vert x - y \Vert^{\lambda}} \; dx.
\end{equation}
\end{definition}
Consequently, for these functions, we have the following property (see Aubin \cite{aubin2012nonlinear}).
\begin{lemma} \label{L2}
Let $u \in C^1_c(\Omega)$. If $\lambda= N-1$, then
\begin{equation*}
\vert u(x) \vert \leq C H(\nabla _1 u)(x),
\end{equation*}
for $x \in \Omega$.
\end{lemma}
Next, we enunciate the following lemma (see Hebey-Vaugon \cite{HEBEY1997859}).
\begin{lemma} [Hebey-Vaugon, 1997] \label{L1}
Let $\Omega$ be a bounded manifold of class $C^1$ in $\mathbb{R}^N$, and let G be a compact subgroup of $\mathcal{O}(N)$. Let $x \in \Omega$ and set $\delta = dim O^{x}_{G}$ with $\delta \geq 1$. There exists a coordinate chart $(A, \phi)$ of $\Omega$ in $x$, with $\phi: \mathbb{R}^N \to \mathbb{R}^N$ a diffeomorphism of class $C^1$ such that:
\begin{enumerate}
\item $\phi(A)= U \times V$, where $U$ is a open subset of $\mathbb{R}^{\delta}$ and $V$ is a open subset of $\mathbb{R}^{N-\delta}$,
\item For all $y \in A$, $U \times \Pi_2 (\phi (y)) \subset \phi(O^{y}_{G} \cap A)$, where $\Pi_2 : \mathbb{R}^{\delta} \times \mathbb{R}^{N-\delta} \to \mathbb{R}^{N-\delta}$ is the projection on the second coordinate.
\end{enumerate}
\end{lemma}
The purpose of this section is to prove the next result:
\begin{lemma} \label{LemaDesiguRi}
Let  $G$ a compact subgroup of $\mathcal{O}(N)$, $\Omega$ a bounded G- invariant open set of class $C^1$ in $\rn$, with $\mbox{dim} \; \Omega =N$; and $k =\displaystyle \min_{x \in \overline{\Omega}} \mbox{dim} \; O^{x}_{G}$. Then there exists a positive constant $C$, such that for all $u \in L^{p}_{G}(\Omega)$, $p>1$ and $1<r \leq \frac{p(N-k)}{N-k-p}$,
$$\Vert H(u) \Vert_{L^r(\Omega)} \leq C \Vert u \Vert_{L^{p}_{G}(\Omega)}$$ where $V$ is the function given in (\ref{V function}).
\end{lemma}
\begin{proof}
Since $\Omega$ is bounded, then $\co$ is compact. Thus  from Lemma \ref{L1}, there exists $\{A_l: l=1, 2, \ldots, \eta \}$ a collection of open subsets such that $\co \subset \displaystyle \bigcup_{l=1}^{\eta} A_l$ and diffeomorphism $\phi_l: A_{l} \to U_l \times V_l \subset \r^{\delta} \times \r^{N - \delta}$ of class $C^1$ for each $l=1, 2, \ldots , \eta$ such that $\phi(\o)$ is covered by a finite number of sets $\phi_l(A_l)$ and holds
\begin{enumerate}
\item $\phi_l(A_l) = U_l \times V_l$, where $U_l$ is a open subset $G$-invariant of $\mathbb{R}^{\delta_l}$ and $V_l$ is an open subset $G$-invariant of $\mathbb{R}^{N-\delta_l}$, with $\delta_l \in \mathbb{N}$ and $\delta_l \geq k$,
\item $U_l$ and $V_l$ are bounded and they have smooth boundary.
\item For all $y \in A_l$, $U_l \times \Pi_2 (\phi_l (y)) \subset \phi_l(O^{x}_{G} \cap A_l)$, where $\Pi_2 : \mathbb{R}^{\delta_l} \times \mathbb{R}^{N-\delta_l} \to \mathbb{R}^{N-\delta_l}$ is the projection on the second coordinate.
\end{enumerate}
We consider $\overline{A} _l \subset \o$, otherwise, we take $A_l \cap \o$, which it is compact because $\Omega$ is bounded. Without loss of generality, we assume that $\phi_l$ is defined on an open and bounded set $\tilde{A}_l$ of class $C^1$, with $\overline{A}_l \subset \tilde{A}_l$, and such that $\phi_l(\tilde{A_l})=\tilde{U}_l \times \tilde{V}_l$ with $\overline{U}_l \subset \tilde{U}_{l}$ y $\overline{V}_l \subset \tilde{V}_{l}$, where $\tilde{U}_{l}$ and $\tilde{V}_{l}$ bounded open sets of class $C^1$.\\
\\Of the above, since $\overline{A}_l$ is compact and $\phi_l$ is a diffeomorphism of class $C^1$, $\phi_l(\overline{A}_l)$ is compact and $\mbox{det} (D \phi_l^{-1}(x))$ is continuous on $\phi_l(\tilde{A}_l)$ consequently has a maximum and minimum, that is, for each $l=1, 2, \ldots , \eta$, there exists positive real numbers $C_l$ and $D_l$ such that
\begin{equation}
C_l \leq \vert \mbox{det} (D \phi_l^{-1}(x)) \vert \leq D_l.
\end{equation}
Let $u \in L^{p}_{G}(\Omega)$. According to item 3, and since $u$ is $G-$invariant, one has that for any $l$, any $x_1, x_{1}' \in U_l$, and $x_2 \in V_l$, $u(\phi_l^{-1}(x_1,x_2))=u(\phi_l^{-1}(x_{1}',x_2))$. Consequently, for $l$ there exists $\tilde{u}_l \in L^{p}(\mathbb{R}^{N-\delta_l})$ such that for $x \in U_l$ and any $y \in V_l$,
\begin{equation}
u(\phi_l^{-1}(x,y)) = \tilde{u}_l(y).
\end{equation}
We have that for a real number $p > 1$,
\begin{equation*}
\begin{split}
\displaystyle \int_{A_l} \vert u(x) \vert ^p \, dx & = \displaystyle \int_{U_l \times V_l} \vert u(\phi_{l}^{-1}(v_1,v_2)) \vert^p  \vert \mbox{det} (D\phi_l^{-1})(v_1, v_2) \vert \: dv_1 \; dv_2\\
& \leq D_l \displaystyle \int_{U_l \times V_l} \vert u(\phi_l^{-1}(v_1,v_2)) \vert^p \; dv_1 \; dv_2 \\
& = D_l \displaystyle \int_{U_l} dv_1  \displaystyle \int_{V_l} \vert \tilde{u}_l(v_2) \vert^p \; dv_2\\
& = D_l \vert U_l \vert \displaystyle \int_{V_l} \vert \tilde{u}_l(v_2) \vert^p \; dv_2\\
& = \tilde{D}_l \displaystyle \int_{V_l} \vert \tilde{u}_l (v_2) \vert^p \; dv_2,
\end{split}
\end{equation*}
where $D_l$ and $\tilde{D}_l$ are positive constants which do not depend of $u$. Thus,
\begin{equation}
\left( \displaystyle \int_{A_l} \vert u(x) \vert ^p \, dx \right)^{\frac{1}{p}} \leq E_l \left( \displaystyle \int_{V_l} \vert \tilde{u}_l (v_2) \vert^p \; dv_2 \right)^{\frac{1}{p}}.
\end{equation}
Similarly, one has that for $l$ and any $p > 1$,
\begin{equation} \label{Porabajo}
\left( \displaystyle \int_{A_l} \vert u(x) \vert ^p \, dx \right)^{\frac{1}{p}} \geq C_l \left( \displaystyle \int_{V_l} \vert \tilde{u}_l (v_2) \vert^p \; dv_2 \right)^{\frac{1}{p}}.
\end{equation}
where $C_l$ is a positive constant which does not depend of $u$.\\
\\Now we choose $p >1$ and $r>1$ such that $\frac{1}{r}=\frac{N-1}{N}+\frac{1}{p}-1$ and denote  $h(y)=H(u)(y)$. From Lemma  \ref{LemaSobolev}, taking $\lambda= N-1$ and $u \in L_{G}^{p}(\Omega)$, for all $g \in L^{q}_{G}(\Omega)$ with $q \geq 1$ and $\frac{1}{r}+\frac{1}{q}=1$, then there exists a constant $K(p,q,N)>0$ such that 
\begin{equation} \label{Sobineq}
\displaystyle \int_{\Omega} h(y) g(y) \, dy \leq \displaystyle \int_{\Omega} \int_{\Omega} \frac{u(x) g(y)}{\Vert x-y \Vert ^ \lambda } \, dxdy \leq  K(p,q,N) \Vert u \Vert_{L^p_{G}(\Omega)} \Vert g \Vert_{L^q_{G}(\Omega)}.
\end{equation}
Hence from the inequality \eqref{Sobineq}, we have
\begin{equation} \label{Sobineq2}
\begin{split}
\Vert h \Vert^{r}_{L^{r}(\o)} & = \displaystyle \io [h(y)]^{r} \, dy\\
& = \displaystyle \io h(y) [h(y)]^{r-1}\\
& \leq K(p,q,N) \Vert u \Vert_{L^p_{G}(\Omega)} \Vert h^{r-1} \Vert_{L^q_{G}(\Omega)},
\end{split}
\end{equation}
since $\frac{1}{r} + \frac{1}{q}=1$, we observe that
\begin{equation} \label{Lq}
\begin{split}
\Vert h^{r-1} \Vert_{L^q_{G}(\Omega)} & = \left(\displaystyle \io  ([h(y)]^{r-1})^q \right)^{\frac{1}{q}}\\
& = \left(\displaystyle \io  ([h(y)]^{r-1})^{\frac{r}{r-1}} \right)^{\frac{r-1}{r}}\\
& = \left(\displaystyle \io  [h(y)]^{r} \right)^{\frac{r-1}{r}}\\
& = \Vert h \Vert_{L^{r}(\o)}^{r-1}.
\end{split}
\end{equation}
Therefore, from \eqref{Sobineq2} and \eqref{Lq},
\begin{equation*}
\begin{split}
\Vert h \Vert^{r}_{L^{r}(\o)} & \leq K(p,q,N) \Vert u \Vert_{L^p_{G}(\Omega)} \Vert h \Vert_{L^{r}(\o)}^{r-1}\\
\Vert h \Vert_{L^{r}(\o)} & \leq K(p,q,N) \Vert u \Vert_{L^p_{G}(\Omega)},
\end{split}
\end{equation*}
thus, $h \in L^{r}(\o)$.\\
\\As above, using a similar argument, since $h \in L_{G}^{r} (\Omega)$, we have that for any $u_1, u_1' \in U_l$, and $u_2 \in V_l$, $h(\phi_l^{-1}(u_1,u_2))=h(\phi_l^{-1}(u_1',u_2))$. Consequently, we can define 
\begin{equation}
\tilde{h}_l(u_2)=h(\phi_l^{-1}(u_1,u_2)).
\end{equation}
such that $\tilde{h}_l \in L^{r}(\mathbb{R}^{N-\delta_l})$.

We have for any real number $r \geq 1$,
\begin{equation*}
\begin{split}
\displaystyle \int_{A_l} \vert h(y) \vert ^r \, dy & = \displaystyle \int_{U_l \times V_l} \vert h(\phi_{l}^{-1}(u_1,u_2)) \vert^r  \vert \mbox{det} (\phi_l^{-1})'(u_1, u_2) \vert \: du_1 \; du_2\\
& \leq D_l \displaystyle \int_{U_l \times V_l} \vert h(\phi_l^{-1}(u_1,u_2)) \vert^r \; du_1 \; du_2 \\
& = D_l \displaystyle \int_{U_l} du_1  \displaystyle \int_{V_l} \vert \tilde{h}_l(u_2) \vert^r \; du_2\\
& = D_l \vert U_l \vert \displaystyle \int_{V_l} \vert \tilde{h}_l(u_2) \vert^r \; du_2\\
& = \tilde{D}_l \displaystyle \int_{V_l} \vert \tilde{h}_l (u_2) \vert^r \; du_2,
\end{split}
\end{equation*}
where $D_l$ and $\tilde{D}_l$ are positive constants which do not depend of $h$. Thus
\begin{equation} \label{Enc1}
\left(\displaystyle \int_{A_l} \vert h(y) \vert ^r \, dy \right)^{\frac{1}{r}} \leq F_l \left( \displaystyle \int_{V_l} \vert \tilde{h}_l (u_2) \vert^r \; du_2 \right)^{\frac{1}{r}}.
\end{equation}
Notice that the first integral is over $A_{l}$ as an open subset from the cover of $\overline{A} \subset \mathbb{R}^N$, while the second one is over $V_{l}$, an open subset of $\mathbb{R}^{N-\delta _l }$.\\
\\ For $\tilde{h}_l \in L^{r}(\mathbb{R}^{N-\delta_l})$, $p >1$ and $r>1$ given above, from lemma \ref{LemaSobolev}, using $\lambda= N- \delta_l -1$, 
for all $\tilde{g} \in L^{q}(\mathbb{R}^{N-\delta_l})$ with $q \geq 1$ and $\frac{1}{r}+\frac{1}{q}=1$, there exists a constant  $K(p,q,N-\delta_l)>0$ such that 
\begin{equation} \label{Sobineq3}
\displaystyle \int_{\mathbb{R}^{N-\delta_l}} \tilde{h}_l(y) \tilde{g}(y) \, dy \leq K(p,q,N-\delta_l) \Vert \tilde{u}_l \Vert_{L^p(\mathbb{R}^{N-\delta_l})} \Vert \tilde{g} \Vert_{L^q(\mathbb{R}^{N-\delta_l})}.
\end{equation}
consequently we have
\begin{equation} \label{Sobineq4}
\begin{split}
\Vert \tilde{h}_l \Vert^{r}_{L^{r}(\r^{N-\delta_{l}})} & = \displaystyle \int_{\r^{N-\delta_l}} [\tilde{h}_l(y)]^{r} \, dy\\
& = \displaystyle \io \tilde{h}_l(y) [\tilde{h}_l(y)]^{r-1}\\
& \leq K(p,q,N-\delta _l) \Vert \tilde{u}_l \Vert_{L^{p}(\r^{N-\delta_{l}})} \Vert \tilde{h}_l^{r-1} \Vert_{L^{q}(\r^{N-\delta_{l}})}.
\end{split}
\end{equation}
Since $\frac{1}{r} + \frac{1}{q}=1$, we observe that
\begin{equation} \label{Lq1}
\begin{split}
\Vert \tilde{h}_l^{r-1} \Vert_{L^{q}(\r^{N-\delta_{l}})} & = \left(\displaystyle \int_{\r^{N-\delta_l}}  ([\tilde{h}_l(y)]^{r-1})^q \right)^{\frac{1}{q}}\\
& = \left(\displaystyle \int_{\r^{N-\delta_l}}  ([\tilde{h}_l(y)]^{r-1})^{\frac{r}{r-1}} \right)^{\frac{r-1}{r}}\\
& = \left(\displaystyle \int_{\r^{N-\delta_l}}  [\tilde{h}_l(y)]^{r} \right)^{\frac{r-1}{r}}\\
& = \Vert \tilde{h}_l \Vert_{L^{r}(\r^{N-\delta_{l}})}^{r-1}.
\end{split}
\end{equation}
Therefore, from \eqref{Sobineq4} and \eqref{Lq1},
\begin{equation*}
\Vert \tilde{h}_l \Vert^{r}_{L^{r}(\r^{N-\delta_{l}})} \leq K(p,q,N-\delta_{l}) \Vert \tilde{u}_l \Vert_{L^{p}(\r^{N-\delta_{l}})} \Vert \tilde{h}_l \Vert_{L^{r}(\r^{N-\delta_{l}})}^{r-1},
\end{equation*}
then
\begin{equation} \label{Desigu001}
\Vert \tilde{h}_l \Vert_{L^r(\mathbb{R}^{N-\delta_l})} \leq K(p,q,N-\delta_l) \Vert \tilde{u}_l \Vert_{L^p(\mathbb{R}^{N-\delta_l})}.
\end{equation}
From the inequalities \eqref{Enc1} and \eqref{Desigu001} we have that
\begin{equation*}
\begin{split}
\left( \displaystyle \int_{A_l} \vert h(y) \vert ^r \, dy \right)^{\frac{1}{r}} & \leq F_l \left( \displaystyle \int_{V_l} \vert \tilde{h}_l (u_2) \vert^r \; du_2 \right)^{\frac{1}{r}} \\
& \leq F_l K \left( \displaystyle \int_{V_l} \vert \tilde{u}_l (v_2) \vert^p \; dv_2 \right)^{\frac{1}{p}}.
\end{split}
\end{equation*}
If $N- \delta_l >p$, for $\frac{1}{r}=\frac{N-\delta_l-1}{N-\delta_l}+\frac{1}{p}-1$, then $r=\frac{(N-\delta_l)p}{N-\delta_l-p}$, for this $r$ and $u \in L^{p}_{G}(A_l)$, follows from \eqref{Porabajo}
\begin{equation} \label{Enc2}
\begin{split}
\left( \displaystyle \int_{A_l} \vert h(y) \vert ^r \, dy \right)^{\frac{1}{r}} & \leq  F_l K \left( \displaystyle \int_{V_l} \vert \tilde{u}_l (v_2) \vert^p \; dv_2 \right)^{\frac{1}{p}}\\
& \leq \frac{F_l}{C_l} K \left( \displaystyle \int_{A_l} \vert u(x) \vert ^p \, dx \right)^{\frac{1}{p}}.
\end{split}
\end{equation}
To extend the rank of $r$ to $1< r \leq \frac{(N - \delta_l)p}{N - \delta_l -p}$, we use that $A_l$ is bounded and the inequality 
\begin{equation} \label{Enc3}
\Vert h \Vert_{L^{r}_{G}(A_l)} \leq \vert A_l \vert^{\frac{r-1}{r}} \Vert h \Vert_{L^{\frac{(N - \delta_l)p}{N - \delta_l -p}}_{G}(A_l)},
\end{equation}
to obtain that \eqref{Enc2} maintains.

Since $N - \delta_l \leq N-k$, then $p^{*}(N, \delta_l, p) \geq p^{*}(N,k,p)$. Moreover, since $\{A_l: l=1, 2, \ldots, \eta \}$ are covering $\overline{\Omega}$, if $X= \displaystyle \bigcup_{l=1}^{\eta} A_l$, then
\begin{equation} \label{Desigu002}
\begin{split}
\left( \displaystyle \int_{\Omega} \vert h(y) \vert^r \, dy \right)^{\frac{1}{r}} & \leq \left( \displaystyle \int_{X} \vert h(y) \vert ^r \, dy \right)^{\frac{1}{r}}\\
& \leq \displaystyle \sum_{l=1}^{\eta} \left(\displaystyle \int_{A_l} \vert h(y) \vert ^r \, dy \right)^{\frac{1}{r}}\\
& \leq \displaystyle \sum_{l=1}^{\eta} \frac{F_l}{C_l} K \left(  \displaystyle \int_{A_l} \vert u(x) \vert ^p \, dx \right)^{\frac{1}{p}}\\
&= \displaystyle \sum_{l=1}^{\eta} K_{l} \left(  \displaystyle \int_{A_l} \vert u(x) \vert ^p \, dx \right)^{\frac{1}{p}}.
\end{split}
\end{equation}
In addition,
\begin{equation} \label{Desigu003}
\begin{split}
K_{l}(p,q,N-\delta_l) \left(\displaystyle \int_{A_l} \vert u(x) \vert ^p \, dx \right)^{\frac{1}{p}} & \leq K_{l}(p,q,N-\delta_l) \left(\displaystyle \int_{\Omega} \vert u(x) \vert ^p \, dx \right)^{\frac{1}{p}}\\
\displaystyle \sum_{l=1}^{\eta} K_{l}(p,q,N-\delta_l) \left(\displaystyle \int_{A_l} \vert u(x) \vert ^p \, dx \right)^{\frac{1}{p}} & \leq \displaystyle \sum_{l=1}^{\eta}  K_{l}(p,q,N-\delta_l) \left(\displaystyle \int_{\Omega} \vert u(x) \vert ^p \, dx \right)^{\frac{1}{p}}\\
\displaystyle \sum_{l=1}^{\eta} K_{l}(p,q,N-\delta_l) \left(\displaystyle \int_{A_l} \vert u(x) \vert ^p \, dx \right)^{\frac{1}{p}} & \leq  A \left(\displaystyle \int_{\Omega} \vert u(x) \vert ^p \, dx \right)^{\frac{1}{p}}.
\end{split}
\end{equation}
Then, from the inequalities \eqref{Desigu002} and \eqref{Desigu003} we have that 
\begin{equation*}
\begin{split}
\left( \displaystyle \int_{\Omega} \vert h(y) \vert^r \, dy \right)^{\frac{1}{r}} & \leq \displaystyle \sum_{l=1}^{\eta} K_l(p,q,m) \left(  \displaystyle \int_{A_l} \vert u(x) \vert ^p \, dx \right)^{\frac{1}{p}}\\
& \leq C \left(\displaystyle \int_{\Omega} \vert u(x) \vert ^p \, dx \right)^{\frac{1}{p}},
\end{split}
\end{equation*}
so that,
\begin{equation}
\Vert V(u) \Vert_{L^r(\Omega)} \leq C \Vert u \Vert_{L^{p}_{G}(\Omega)}.
\end{equation}
\end{proof}

\section{Proofs of main results} \label{sec4}

In this section we give the proofs of the theorem \ref{teo1} and theorem \ref{teo2} which we discussed above.
\begin{proof}[Proof of Theorem \ref{teo1}]
Let $0<a<2$ and $0<b<2$ be real numbers. For $u \in H^1_{s,G}(\Omega)$, from lemmas \ref{Lema1W} and \ref{L2},
\begin{equation*}
\begin{split}
\vert u(x) \vert^b & \leq \dfrac{C_j^b}{\Vert x_2 \Vert^{\frac{b(j-2)}{2}}} \left( \displaystyle \int_{\Omega_2} \vert \nabla  u \vert^2 \right)^{\frac{b}{2}}\\
\vert u(x) \vert^a & \leq C^a (H(\nabla_1 u))^a,\\
\vert u(x) \vert^{a+b} & \leq \dfrac{C^a C_j^b (H(\nabla_1 u))^a}{\Vert x_2 \Vert^{\frac{b(j-2)}{2}}} \left( \displaystyle \int_{\Omega_2} \vert \nabla  u \vert^2 \right)^{\frac{b}{2}}\\
\Vert x_2 \Vert^l \vert u(x) \vert^{a+b} & \leq \dfrac{C^a C_j^b (H(\nabla_1 u))^a}{\Vert x_2 \Vert^{\frac{b(j-2)}{2}-l}} \left( \displaystyle \int_{\Omega_2} \vert \nabla  u \vert^2 \right)^{\frac{b}{2}}.
\end{split}
\end{equation*}
We set $\gamma=\frac{b(j-2)}{2}-l$. Then,
\begin{equation} \label{Desigu1}
\Vert x_2 \Vert^l \vert u(x) \vert^{a+b} \leq \dfrac{B (H(\nabla_1 u))^a}{\Vert x_2 \Vert^{\gamma}} \left( \displaystyle \int_{\Omega_2} \vert \nabla  u \vert^2 \right)^{\frac{b}{2}}.
\end{equation}
Integrating \eqref{Desigu1} on $\Omega_1$
\begin{equation*}
\displaystyle \int_{\Omega_1} \Vert x_2 \Vert^l \vert u(x) \vert^{a+b} \leq \frac{B}{\Vert x_2 \Vert^{\gamma}} \displaystyle \int_{\Omega_1} \left[ (H(\nabla_1 u))^a \left( \displaystyle \int_{\Omega_2} \vert \nabla  u \vert^2 \right)^{\frac{b}{2}} \right].
\end{equation*}
Applying Hölder’s inequality with conjugated exponents $\frac{2}{b}$ and $\frac{2}{2-b}$ we have
\begin{equation} \label{Holder1}
\displaystyle \int_{\Omega_1} \Vert x_2 \Vert^l \vert u(x) \vert^{a+b} \leq \frac{B}{\Vert x_2 \Vert^{\gamma}} \left( \displaystyle \int_{\Omega_1}  (H(\nabla_1 u))^{\frac{2a}{2-b}} \right)^{\frac{2-b}{2}} \left( \displaystyle \int_{\Omega} \vert \nabla  u \vert^2  \right)^{\frac{b}{2}}.
\end{equation}
Rewriting 
\begin{equation*}
\begin{split}
\left( \displaystyle \int_{\Omega_1}  (H(\nabla_1 u))^{\frac{2a}{2-b}} \right)^{\frac{2-b}{2}} & = \left( \left( \displaystyle \int_{\Omega_1}  (H(\nabla_1 u))^{\frac{2a}{2-b}} \right)^{\frac{2-b}{2a}} \right)^a\\
& = \Vert H(\nabla_1 u) \Vert^{a}_{L^{\frac{2a}{2-b}}(\Omega_1)}
\end{split}
\end{equation*}
and applying lemma \ref{LemaDesiguRi}, setting $u$ as $\nabla_1 u$ and $p=2$, we have
\begin{equation*}
\left( \displaystyle \int_{\Omega_1}  (H(\nabla_1 u))^{\frac{2a}{2-b}} \right)^{\frac{2-b}{2}} = \Vert H(\nabla_1 u) \Vert^{a}_{L^{\frac{2a}{2-b}}(\Omega_1)} \leq A \left( \displaystyle \int_{\Omega_1} \vert \nabla_1 u \vert^2 \right)^{\frac{a}{2}}.
\end{equation*}
Notice from lemma \ref{LemaDesiguRi}, must be satisfied $\frac{2a}{2-b} \leq \frac{2(m-k)}{m-k-2}$, or equivalently
\begin{equation} \label{Condi1}
(m-k-2)a+(m-k)b \leq 2(m-k).
\end{equation}
Under this condition, from \eqref{Holder1}, we have
\begin{equation}
\displaystyle \int_{\Omega_1} \Vert x_2 \Vert^l \vert u(x) \vert^{a+b} \leq \frac{B}{\Vert x_2 \Vert^{\gamma}} A \left( \displaystyle \int_{\Omega_1} \vert \nabla_1 u \vert^2 \right)^{\frac{a}{2}} \left( \displaystyle \int_{\Omega} \vert \nabla  u \vert^2  \right)^{\frac{b}{2}}.
\end{equation}
Now integrating on $\Omega_2$
\begin{equation*}
\begin{split}
\displaystyle \int_{\Omega_2} \displaystyle \int_{\Omega_1} \Vert x_2 \Vert^l \vert u(x) \vert^{a+b} & \leq D \displaystyle \int_{\Omega_2} \left[ \Vert x_2 \Vert^{-\gamma} \left( \displaystyle \int_{\Omega_1} \vert \nabla_1 u \vert^2 \right)^{\frac{a}{2}} \left( \displaystyle \int_{\Omega} \vert \nabla  u \vert^2  \right)^{\frac{b}{2}} \right]\\
& = D \left( \displaystyle \int_{\Omega} \vert \nabla  u \vert^2  \right)^{\frac{b}{2}} \displaystyle \int_{\Omega_2} \left[ \Vert x_2 \Vert^{-\gamma} \left( \displaystyle \int_{\Omega_1} \vert \nabla_1 u \vert^2 \right)^{\frac{a}{2}} \right],
\end{split}
\end{equation*}
then
\begin{equation*}
\displaystyle \int_{\Omega} \Vert x_2 \Vert^l \vert u(x) \vert^{a+b} \leq D \left( \displaystyle \int_{\Omega} \vert \nabla  u \vert^2  \right)^{\frac{b}{2}} \displaystyle \int_{\Omega_2} \left[ \Vert x_2 \Vert^{-\gamma} \left( \displaystyle \int_{\Omega_1} \vert \nabla_1 u \vert^2 \right)^{\frac{a}{2}} \right].
\end{equation*}
Again by Hölder's inequality with conjugated exponents $\frac{2}{a}$ and $\frac{2}{2-a}$ we have
\begin{equation} \label{Holder2}
\displaystyle \int_{\Omega} \Vert x_2 \Vert^l \vert u(x) \vert^{a+b} \leq D \left(\displaystyle \int_{\Omega_2} \Vert x_2 \Vert^{\frac{-2 \gamma}{2-a}} \right)^{\frac{2-a}{2}} \left( \displaystyle \int_{\Omega} \vert \nabla  u \vert^2  \right)^{\frac{a+b}{2}}.
\end{equation}
We will look the behavior of the integral $\left(\displaystyle \int_{\Omega_2} \Vert x_2 \Vert^{\frac{-2 \gamma}{2-a}} \right)^{\frac{2-a}{2}}$. We consider the equalities
\begin{equation} \label{Integralfinita1}
\begin{split}
\displaystyle \int_{\Omega_2} \Vert x_2 \Vert^{\frac{-2 \gamma}{2-a}} & = j \omega_j \displaystyle \int_{0}^{R} t^{\frac{-2 \gamma}{2-a}+j-1} \, dt\\
& = j \omega_j \dfrac{1}{\frac{-2 \gamma}{2-a}+j} R^{\frac{-2 \gamma}{2-a}+j};
\end{split}
\end{equation}
then this integral is finite if and only if $\frac{-2 \gamma}{2-a}+j>0$, that is, $\frac{2 \gamma}{2-a}<j$, or equivalently,
\begin{equation} \label{Condi2}
ja+(j-2)b<2(j+l),
\end{equation}
then the inequality \eqref{Holder2} is reduced to
\begin{equation} \label{Desigu3}
\displaystyle \int_{\Omega} \Vert x_2 \Vert^l \vert u(x) \vert^{a+b} \leq A \left( \displaystyle \int_{\Omega} \vert \nabla  u \vert^2  \right)^{\frac{a+b}{2}}.
\end{equation}
We can see that the inequality \eqref{Desigu3} is satisfied if and only if the conditions \eqref{Condi1} and \eqref{Condi2} are satisfied. Both conditions are equivalents to require that $0<a<2$ and $0<b<2$ solve the following inequalities  
\begin{equation} \label{Desigu4}
\left\lbrace \begin{matrix}
(m-k-2)a+(m-k)b \leq 2(m-k)\\
ja+(j-2)b <2(j+l),
\end{matrix} \right.
\end{equation}
or equivalently,
\begin{equation}
0<a \leq \frac{(m-k)(l+2)}{(m-k)+j-2}, \ \ \; 0<b < \frac{2(m-k)-(m-k-2)a}{m-k},
\end{equation} 
\begin{equation} \label{Puntodeinterseccion}
\frac{(m-k)(l+2)}{m-k+j-2} \leq a <2, \ \ \; 0<b< \frac{2(j+l)-ja}{j-2}.
\end{equation}
From inequality \eqref{Puntodeinterseccion}, since $\frac{(m-k)(l+2)}{m-k+j-2} <2$, then $l<\frac{2j-4}{m-k}$. We have to choose the appropriate conditions such that $a$ and $b$ satisfy \eqref{Desigu4}, such that $a+b$ approaches $2^{*}_{k} + \psi$. We define
\begin{equation*}
l^{*}= \min \left\lbrace l, \frac{2j-4}{m-k} \right\rbrace,
\end{equation*}
\begin{equation*}
\psi= \frac{2 l^{*}}{m-k+j-2}.
\end{equation*}
If $m \geq 3$, we take $l' \in (0, l^{*})$, such that
$$a=\frac{(m-k)(l'+2)}{m-k+j-2},$$
$$b= \frac{2(j+l')-l'(m-k)}{m-k+j-2},$$
then $a$ and $b$ satisfy \eqref{Desigu4}. Moreover
\begin{equation*}
\begin{split}
a+b & = \frac{(m-k)(l'+2)}{m-k+j-2} + \frac{2(j+l')-l'(m-k)}{m-k+j-2}\\
& = \frac{2(m-k)+2j+2l'}{m-k+j-2}\\
& = \frac{2(m-k+j)}{m-k+j-2}+ \frac{2l'}{m-k+j-2}\\
& = \frac{2(N-k)}{N-k-2}+ \frac{2l'}{N-k-2}.
\end{split}
\end{equation*}
If $l' \to l^{*}$, then $a+b \to 2^{*}_{k} + \psi$.\\
\\If $m=2$, since $1 \leq k<m$, then $k=1$. We take $l' \in (0, l^{*})$ as in the previous case, we have that $\frac{2+l'}{j-1}<2$. Moreover, we have
$$a=\frac{l'+2}{j-1},$$
$$b=\frac{2j+l'}{j-1},$$
and
\begin{equation*}
\begin{split}
a+b & = \frac{l'+2}{j-1} + \frac{2j+l'}{j-1}\\
& = \frac{l'+2+2j+l'}{j-1}\\
& = \frac{2+2j+2l'}{j-1}\\
& = \frac{2j+2}{j-1}+\frac{2l'}{j-1}\\
& = \frac{2(N-k)}{N-k-2}+ \frac{2l'}{j-1}.
\end{split}
\end{equation*}
If $l' \to l^{*}$, then $a+b \to 2^{*}_k + \psi$, that is, $a+b\in (1,2^{*}_k+\psi )$. Making $a+b=q$ in  (\ref{Desigu3})  we can see that the embedding $H^1_{s,G}(\Omega) \subset L^q_h(\Omega)$ is continuous.\\
\\ Finally, to conclude that the embedding is compact, we can follow a similar procedure for the proof of Kondrachov's theorem, which can be found in \cite{brezis2011functional}. This same argument is also applied to theorem \ref{teo2} of this paper.
\end{proof}
\begin{proof} [Proof of Theorem \ref{teo2}]
Let $\alpha, \beta \in (0,1)$ and $b>1$ be real numbers. For $u \in H^1_{s,G}(\Omega) \cap C^1_{c,G}(\Omega)$, from lemmas \ref{Lema2W} and \ref{L2},
\begin{equation*}
\begin{split}
(\vert u(x) \vert^b)^{\beta} & \leq \dfrac{(b A_k)^{\beta}}{\Vert x_2 \Vert^{\beta}} \left( \displaystyle \int_{\Omega_2} \vert u \vert^{b-1} \vert \nabla  u \vert \right)^{\beta}\\
(\vert u(x) \vert^b)^{\alpha} & \leq C^{\alpha} (H(\vert u \vert^{b-1} \nabla_1 u))^{\alpha},\\
\vert u(x) \vert^{b(\alpha + \beta)} & \leq \dfrac{C^{\alpha} (bA_k)^{\beta} (H(\vert u \vert^{b-1} \nabla_1 u))^{\alpha}}{\Vert x_2 \Vert^{\beta}} \left( \displaystyle \int_{\Omega_2} \vert u \vert^{b-1} \vert \nabla  u \vert \right)^{\beta}\\
\Vert x_2 \Vert^l \vert u(x) \vert^{b(\alpha + \beta)} & \leq \dfrac{E (H(\vert u \vert^{b-1} \nabla_1 u))^{\alpha}}{\Vert x_2 \Vert^{\beta - l}} \left( \displaystyle \int_{\Omega_2} \vert u \vert^{b-1} \vert \nabla  u \vert \right)^{\beta}.
\end{split}
\end{equation*}
Integrating on $\Omega_1$ we have that
\begin{equation*}
\displaystyle \int_{\Omega_1} \Vert x_2 \Vert^l \vert u(x) \vert^{b(\alpha + \beta)} \leq \frac{E}{\Vert x_2 \Vert^{\beta -l}} \displaystyle \int_{\Omega_1} \left[(H(\vert u \vert^{b-1} \nabla_1 u))^{\alpha} \left( \displaystyle \int_{\Omega_2} \vert u \vert^{b-1} \vert \nabla  u \vert \right)^{\beta} \right].
\end{equation*}
Since $\beta <1$, $\frac{1}{\beta}>1$. Applying the Hölder's inequality with conjugated exponents $\frac{1}{\beta}$ y $\frac{1}{1-\beta} $ we have
\begin{equation} \label{Holder3}
\displaystyle \int_{\Omega_1} \Vert x_2 \Vert^l \vert u(x) \vert^{b(\alpha + \beta)} \leq \frac{E}{\Vert x_2 \Vert^{\beta-l}} \left( \displaystyle \int_{\Omega_1}  (H(\vert u \vert^{b-1} \nabla_1 u))^{\frac{\alpha}{1-\beta}} \right)^{1-\beta} \left( \displaystyle \int_{\Omega} \vert u \vert^{b-1} \vert \nabla  u \vert  \right)^{\beta}.
\end{equation}
We make $\gamma= \beta -l$. Analyzing the first integral on the right side of the inequality \eqref{Holder3}
\begin{equation*}
\begin{split}
\left( \displaystyle \int_{\Omega_1}  (H(\vert u \vert^{b-1} \nabla_1 u))^{\frac{\alpha}{1-\beta}} \right)^{1-\beta} & = \left( \left( \displaystyle \int_{\Omega_1}  (H(\vert u \vert^{b-1} \nabla_1 u))^{\frac{\alpha}{1-\beta}} \right)^{\frac{1-\beta}{\alpha}} \right)^{\alpha}\\
& = \Vert H(\vert u \vert^{b-1} \nabla_1 u) \Vert^{\alpha}_{L^{\frac{\alpha}{1-\beta}}(\Omega_1)}.
\end{split}
\end{equation*}
From lemma \ref{LemaDesiguRi}, with $u= \nabla_1 u$ and $p=1$, we have
\begin{equation*}
\left( \displaystyle \int_{\Omega_1}  (H(\vert u \vert^{b-1} \nabla_1 u))^{\frac{\alpha}{1-\beta}} \right)^{1-\beta} = \Vert H(\vert u \vert^{b-1} \nabla_1 u) \Vert^{\alpha}_{L^{\frac{\alpha}{1-\beta}}(\Omega_1)} \leq A \left( \displaystyle \int_{\Omega_1} \vert u \vert^{b-1} \vert \nabla_1 u \vert \right)^{\alpha},
\end{equation*}
if and only if $\frac{\alpha}{1-\beta} \leq \frac{m-k}{m-k-1}$, or equivalently
\begin{equation} \label{Condi3}
\alpha(m-k-1)+\beta(m-k) \leq m-k.
\end{equation}
From this condition, from \eqref{Holder3}, we have
\begin{equation} \label{Desigu5}
\displaystyle \int_{\Omega_1} \Vert x_2 \Vert^l \vert u(x) \vert^{b(\alpha + \beta)} \leq \frac{D}{\Vert x_2 \Vert^{\gamma}} \left( \displaystyle \int_{\Omega_1} \vert u \vert^{b-1} \vert \nabla_1 u \vert \right)^{\alpha} \left( \displaystyle \int_{\Omega} \vert u \vert^{b-1} \vert \nabla  u \vert  \right)^{\beta}.
\end{equation}
Integrating \eqref{Desigu5} on $\Omega_2$
\begin{equation*}
\begin{split}
\displaystyle \int_{\Omega_2} \displaystyle \int_{\Omega_1} \Vert x_2 \Vert^l \vert u(x) \vert^{b(\alpha + \beta)} & \leq D \displaystyle \int_{\Omega_2} \left[ \Vert x_2 \Vert^{-\gamma} \left( \displaystyle \int_{\Omega_1} \vert u \vert^{b-1} \vert \nabla_1 u \vert \right)^{\alpha} \left( \displaystyle \int_{\Omega} \vert u \vert^{b-1} \vert \nabla  u \vert  \right)^{\beta} \right]\\
& = D \left( \displaystyle \int_{\Omega} \vert u \vert^{b-1} \vert \nabla  u \vert  \right)^{\beta} \displaystyle \int_{\Omega_2} \left[ \Vert x_2 \Vert^{-\gamma} \left( \displaystyle \int_{\Omega_1} \vert u \vert^{b-1} \vert \nabla_1 u \vert \right)^{\alpha} \right],
\end{split}
\end{equation*}
therefore
\begin{equation*}
\displaystyle \int_{\Omega} \Vert x_2 \Vert^l \vert u(x) \vert^{b(\alpha + \beta)} \leq D \left( \displaystyle \int_{\Omega} \vert u \vert^{b-1} \vert \nabla  u \vert  \right)^{\beta} \displaystyle \int_{\Omega_2} \left[ \Vert x_2 \Vert^{-\gamma} \left( \displaystyle \int_{\Omega_1} \vert u \vert^{b-1} \vert \nabla_1 u \vert \right)^{\alpha} \right].
\end{equation*}
Applying the Hölder's inequality using $\frac{1}{\alpha}$ and $\frac{1}{1-\alpha}$ we have
\begin{equation} \label{Holder4}
\displaystyle \int_{\Omega} \Vert x_2 \Vert^l \vert u(x) \vert^{b(\alpha + \beta)} \leq D \left(\displaystyle \int_{\Omega_2} \Vert x_2 \Vert^{\frac{- \gamma}{1-\alpha}} \right)^{1-\alpha} \left( \displaystyle \int_{\Omega} \vert u \vert^{b-1} \vert \nabla u \vert  \right)^{\alpha + \beta}.
\end{equation}
In a similar way than the equalities \eqref{Integralfinita1}, analyzing the behavior of the integral $\left(\displaystyle \int_{\Omega_2} \Vert x_2 \Vert^{\frac{- \gamma}{1-\alpha}} \right)^{1-\alpha}$
\begin{equation*}
\begin{split}
\displaystyle \int_{\Omega_2} \Vert x_2 \Vert^{\frac{- \gamma}{1-\alpha}} & = 2 \omega_2 \displaystyle \int_{0}^{R} t^{\frac{- \gamma}{1-\alpha}+1} \, dt\\
& = 2 \omega_2 \dfrac{1}{\frac{- \gamma}{1-\alpha}+2} R^{\frac{- \gamma}{1-\alpha}+2};
\end{split}
\end{equation*}
we see that this integral is finite if and only if $\frac{- \gamma}{1-\alpha}<2$, or equivalently
\begin{equation} \label{Condi4}
2\alpha+\beta<2+l,
\end{equation}
then \eqref{Holder4} is reduced to
\begin{equation} \label{Desigu6}
\displaystyle \int_{\Omega} \Vert x_2 \Vert^l \vert u(x) \vert^{b(\alpha + \beta)} \leq B \left( \displaystyle \int_{\Omega} \vert u \vert^{b-1} \vert \nabla u \vert  \right)^{\alpha + \beta}.
\end{equation}
Using the Hölder's inequality in \eqref{Desigu6}, we have
\begin{equation*}
\begin{split}
\displaystyle \int_{\Omega} \Vert x_2 \Vert^l \vert u(x) \vert^{b(\alpha + \beta)}  & \leq B \left( \displaystyle \int_{\Omega} \vert u \vert^{b-1} \vert \nabla u \vert  \right)^{\alpha + \beta}\\
& \leq B \left(\displaystyle \int_{\Omega} \vert u \vert^{2(b-1)} \right)^{\frac{\alpha + \beta}{2}} \left(\displaystyle \int_{\Omega} \vert \nabla u \vert^2 \right)^{\frac{\alpha + \beta}{2}}.
\end{split}
\end{equation*}
If $2(b-1) \leq \frac{2N}{N-2}$, then, since $N=m+2$, we have
\begin{equation} \label{Condi5}
b \leq \frac{2(m+1)}{m}.
\end{equation}
From the Sobolev's embedding theorem, we have that if $2(b-1) \leq 2^{*}$, then
\begin{equation*}
\begin{split}
\left(\displaystyle \int_{\Omega} \vert u \vert^{2(b-1)} \right)^{\frac{1}{2(b-1)}} & \leq C \left( \displaystyle \int_{\Omega} \vert \nabla u \vert^2 \right)^{\frac{1}{2}}\\
\left(\displaystyle \int_{\Omega} \vert u \vert^{2(b-1)} \right)^{\frac{1}{(b-1)}} & \leq C  \displaystyle \int_{\Omega} \vert \nabla u \vert^2
\end{split}
\end{equation*}
\begin{equation} \label{Desigu8}
\displaystyle \int_{\Omega} \vert u \vert^{2(b-1)} \leq C^{b-1} \left( \displaystyle \int_{\Omega} \vert \nabla u \vert^2 \right)^{b-1}.
\end{equation}
From the inequalities \eqref{Desigu6}, \eqref{Desigu8} we have
\begin{equation} \label{Desigu9}
\begin{split}
\displaystyle \int_{\Omega} \Vert x_2 \Vert^l \vert u(x) \vert^{b(\alpha + \beta)} & \leq BC^{b-1} \left(\displaystyle \int_{\Omega} \vert \nabla u \vert^2 \right)^{(b-1)\frac{\alpha + \beta}{2}} \left(\displaystyle \int_{\Omega} \vert \nabla u \vert^2 \right)^{\frac{\alpha + \beta}{2}} \\
& =A \left(\displaystyle \int_{\Omega} \vert \nabla u \vert^2 \right)^{\frac{b(\alpha + \beta)}{2}}.
\end{split}
\end{equation}
If the conditions \eqref{Condi3}, \eqref{Condi4} and \eqref{Condi5} are satisfied then we have \eqref{Desigu9}. These conditions are equivalent to require that $0<\alpha<1$, $0<\beta<1$ and $b>1$ satisfy the following
\begin{equation} \label{Desigu12}
\left\lbrace \begin{matrix}
\alpha(m-k-1)+\beta(m-k) \leq m-k\\
2 \alpha+\beta <2+l\\
b \leq \frac{2(m+1)}{m}
\end{matrix} \right.
\end{equation}
or equivalently, that $\alpha$, $\beta$ satisfy
\begin{equation}
0<\alpha \leq \frac{(m-k)(1+l)}{m-k+1}, \; 0<\beta < 1-\frac{m-k-1}{m-k} \alpha,
\end{equation}
\begin{equation} \label{Puntodeinterseccion2}
\frac{(m-k)(1+l)}{m-k+1} \leq \alpha <1, \; 0<\beta< 2+l-2\alpha.
\end{equation}
Since the inequality \eqref{Puntodeinterseccion2} is satisfied, $\frac{(m-k)(1+l)}{m-k+1} <1$, therefore $l<\frac{1}{m-k}$. What remains for us is to choose the appropriate conditions such that $\alpha$ and $\beta$ check \eqref{Desigu12}. Then we define
\begin{equation*}
l^{*}= \min \left\lbrace l, \frac{1}{m-k} \right\rbrace,
\end{equation*}
\begin{equation*}
\psi= -\frac{2k(N-k)}{m(m-k)(m-k+1)} +\frac{2(m+1)}{m(m-k+1)} l^{*},
\end{equation*}
and for $\eta >0$, we make
$$\alpha=\frac{(m-k)(1+l^{*})-\eta}{m-k+1},$$
$$\beta= \frac{2-(m-k-1)l^{*}}{m-k+1},$$
$$b=\frac{2(m+1)}{m},$$
then $\alpha$, $\beta$ and $b$ satisfy \eqref{Desigu12}. 
In addition,
\begin{multline*}
\begin{split}
b(\alpha + \beta) & = \frac{2(m+1)}{m} \left[\frac{(m-k)(1+l^{*})-\eta + 2-(m-k-1)l^{*}}{m-k+1} \right]\\
& = \frac{2(m+1)}{m} \left[\frac{m-k+l^{*}+2-\eta}{m-k+1} \right]\\
& = \frac{2(m+1)}{m} \left[\frac{m-k+2}{m-k+1}+\frac{l^{*}}{m-k+1}-\frac{\eta}{m-k+1} \right]\\
& = \frac{(m+1)}{m} \left[\frac{2(N-k)}{m-k+1}+\frac{2l^{*}}{m-k+1}-\frac{2\eta}{m-k+1} \right]\\
& =\frac{2(N-k)}{N-k-2} - \frac{2(N-k)}{N-k-2} + \frac{2(N-k)(m+1)}{m(m-k+1)}\\
& \hspace{3mm} +\frac{2(m+1)}{m(m-k+1)} l^{*}-\frac{2(m+1)\eta}{m(m-k+1)}\\
& = \frac{2(N-k)}{N-k-2} - \frac{2k(N-k)}{m(m-k)(m-k+1)}\\
& \hspace{3mm} + \frac{2(m+1)}{m(m-k+1)} l^{*}-\frac{2(m+1)\eta}{m(m-k+1)}\\
& = 2^{*}_{k} + \psi-\frac{2(m+1)\eta}{m(m-k+1)}.
\end{split} 
\end{multline*}
If $\eta \to 0$ then $b(\alpha + \beta) \to 2^{*}_{k} + \psi$, that is, $b(\alpha + \beta)\in (1,2^{*}_{k}+\psi)$, making $b(\alpha + \beta)=q$ we can see that the embedding $H^1_{s,G}(\Omega) \subset L^q_h(\Omega)$ is continuous.
\end{proof}

\end{document}